\documentclass[12pt]{article}
\usepackage{amssymb,amsmath}
\begin{document}
\newtheorem{theorem}{Theorem}
\newtheorem{lemma}{Lemma}
\newtheorem{proposition}{Proposition}
\newtheorem{definition}{Definition}
\begin{center}{\Large Introductory Clifford Analysis}
\end{center}
\begin{center}{by John Ryan}\end{center}
\begin{center}{Department of Mathematics, University of Arkansas}
\end{center}
\begin{center}{Fayetteville, AR 72701, U. S. A.}\end{center}
\begin{center}{\Large Introduction}\end{center}

\ We want in this chapter to regard Clifford algebras as natural generalizations of the complex number system. First let us note that if $z$ is a complex number then
$\overline{z}z=\|z\|^{2}$. Now for a quaternion $q$ we also have
$\overline{q}q=\|q\|^{2}$. In this way quaternions may be regarded
as a generalization of the complex number system. It seems natural
to ask if one can extend basic results of one complex variable
analysis on holomorphic function theory to four dimensions using
quaternions. The answer is yes. This was developed by the Swiss
mathematician Rudolph Fueter in the 1930's and 1940's and also my Moisil and Teodorecu \cite{mt}. See for
instance \cite{f}. An excellent review of this work is given in the
survey article "Quaternionic analysis" by A. Sudbery, see \cite{su}. There is 
also earlier work of Dixon \cite{di}.
However in previous lectures we have seen that for a vector $x\in
R^{n}$ when we consider $R^{n}$ as embedded in the Clifford
algebra $Cl_{n}$ then $x^{2}=-\|x\|^{2}$. So it seems reasonable
to ask if all that is known in the quaternionic setting extends to
the Clifford algebra setting. Again the answer is yes. The earlier aspects 
of this study was developed by amongst others Richard
Delanghe \cite{d}, Viorel Iftimie \cite{i} and David Hestenes \cite{h}. The subject that
has grown from these works is now called Clifford analysis.

\ In more recent times Clifford analysis has found a wealth of unexpected
applications in a number of branches of mathematical analysis particularly
classical harmonic analysis, see for instance the work of Alan McIntosh
and his collaborators, for instance \cite{lms,lmq}, Marius Mitrea \cite{m1,m2}
and papers in \cite{r4}. Links to representation theory and several complex variables may be
found in \cite{gm,r1,r2,r3} and elsewhere.

\ The purpose of this paper is to present a review of many of the basic
aspects of Clifford analysis. 

\ Alternative accounts of much of this work together with other related
results can be found in \cite{bds,dss,gs,gm,ks,o,r4,r5}.

\begin{center}{\Large Foundations of Clifford Analysis}\end{center}

\ We start by replacing the vector $x=x_{1}e_{1}+\ldots+x_{n}e_{n}$ by the
differential operator $D=\Sigma_{j=1}^{n}e_{j}\frac{\partial}{\partial x_{j}}$.
One basic but interesting property of $D$ is that $D^{2}=-\triangle_{n}$,
the Laplacian $\Sigma_{j=1}^{n}\frac{\partial^{2}}{\partial x_{j}^{2}}$ in
${\bf{R}}^{n}$. The differential operator $D$ will be called a Dirac operator.
This is because the classical Dirac operator constructed over four dimensional
Minkowski space squares to give the wave operator.

\begin{definition}
Suppose that $U$ is a domain in ${\bf{R}}^{n}$ and $f$ and $g$ are $C^{1}$
functions defined on $U$ and taking values in $Cl_{n}$. Then $f$ is called a
left monogenic function if $Df=0$ on $U$ while $g$ is called a right
monogenic function on $U$ if $gD=0$ where $gD=\Sigma_{j=1}^{n}
\frac{\partial g}{\partial x_{j}}e_{j}$.
\end{definition}

\ Alternatively left monogenic functions are also called left regular
functions and perhaps most appropriately left Clifford holomorphic
functions. The term Clifford holomorphic functions or Clifford analytic functions appears to be due to Stephen Semmes, see \cite{se} and elsewhere. We shall most often use the term Clifford holomorphic functions here.

\ Examples of such functions include the gradients of real valued harmonic
functions on $U$. So if $h$ is harmonic on $U$ and it is also real valued
then $Dh$ is a vector valued left monogenic function. It is also a right
monogenic function. Such a function is more commonly referred to as a
conjugate harmonic function or a harmonic $1$-form. See for instance
\cite{sw}. An example of such a function is $G(x)=\frac{x}{\|x\|^{n}}$.

\ It should be noted that if $f$ and $g$ are left monogenic functions then 
due to the lack of commutativity of the Clifford algebra, it is not in 
general true that their product $f(x)g(x)$ is left monogenic.

\ To introduce other possible examples of left monogenic functions suppose
that $\mu$ is a $Cl_{n}$ valued measure with compact support $[\mu]$ in
${\bf{R}}^{n}$. Then the convolution $\int_{[\mu]}G(x-y)d\mu(y)$ defines a
left monogenic function on the maximal domain lying in
${\bf{R}}^{n}\backslash[\mu]$. The previously defined integral is the 
Cauchy transform of the measure $[\mu]$.

\ Another way to construct examples of left monogenic functions was introduced
by Littlewood and Gay in \cite{lg} for the case $n=3$ and independently
re-introduced for all $n$ by Sommen \cite{s2}. Suppose $U'$ is a domain in
${\bf{R}}^{n-1}$, the span of $e_{2},\ldots,e_{n}$. Suppose also that
$f'(x')$ is a $Cl_{n}$ valued function such that at each point $x'\in U'$
there is a multiple series expansion in $x_{2},\ldots,x_{n}$ that converges
uniformly on some neighbourhood of $x'$ in $U'$ to $f'$. Such a function is
called a real analytic function. The series
\[\Sigma_{k=0}^{\infty}\frac{1}{k!}x_{1}^{k}(-e_{1}D'f'(x')=\exp(-x_{1}e_{1}D')
f'(x')\]
where $D'=\Sigma_{j=2}^{n}e_{j}\frac{\partial}{\partial x_{j}}$, defines a
left monogenic function $f$ in some neighborhood $U(f')$ in ${\bf{R}}^{n}$ of
$U'$. The left monogenic function $f$ is the Cauchy-Kowalewska extension of
$f'$.

\ It should be noted that if $f$ is a left monogenic function then
$\overline{f}$ and $\tilde{f}$ are both right monogenic functions.

\ We now turn to analogues of Cauchy's Theorem and Cauchy's integral formula.
\begin{theorem}{\bf{(The Clifford-Cauchy Theorem):}} Suppose that $f$ is
a
left Clifford holomorphic function on $U$ and $g$ is a right Clifford holomorphic function on $U$.
Suppose also that $V$ is a bounded subdomain of $U$ with piecewise
differentiable boundary $S$ lying in $U$. Then
\begin{equation}
\int_{S}g(x)n(x)f(x)d\sigma(x)=0
\end{equation}
where $n(x)$ is the outward pointing normal vector to $S$ at $x$ and $\sigma$
is the Lebesgue measure on $S$.
\end{theorem}

\ The proof follows directly from Stokes' Theorem. One important point to
keep
in mind though is that as $Cl_{n}$ is not a commutative algebra then it is
important to place the vector $n(x)$ between $f$ and $g$. One then has that
\[\int_{S}g(x)n(x)f(x)d\sigma(x)=\int_{V}((g(x)D)f(x)+g(x)(Df(x)))dx^{n}=0.\]

\ Suppose that $g$ is the gradient of a real valued harmonic function and
$f=1$. Then the real part of Equation 1 gives the following well known integral
formula.
\[\int<grad g(x),n(x)>d\sigma(x)=0.\]

\ We now turn to the analogue of a Cauchy integral formula.
\begin{theorem}{\bf{(Clifford-Cauchy Integral Formula):}}
Suppose that $U$, $V$, $S$, $f$ and $g$ are all as in Theorem 1 and that
$y\in V$. Then
\[f(y)=\frac{1}{\omega_{n}}\int_{S}G(x-y)n(x)f(x)d\sigma(x)\]
and
\[g(y)=\frac{1}{\omega_{n}}\int_{S}g(x)n(x)G(x-y)d\sigma(x)\]
where $\omega_{n}$ is the surface area of the unit sphere in ${\bf{R}}^{n}$.
\end{theorem}
{\bf{Proof:}} The proof follows very similar lines to the argument in one
variable complex analysis. We shall establish the formula for $f(y)$ the proof
being similar for $g(y)$. First let us take a sphere $S^{n-1}(y,r)$ centered
at $y$ and of radius $r$. The radius $r$ is chosen sufficently small so that
the closed disc with boundary $S^{n-1}(y,r)$ lies in $V$. Then by the
Clifford-Cauchy theorem
\[\int_{S}G(x-y)n(x)f(x)d\sigma(x)=
\int_{S^{n-1}(y,r)}G(x-y)n(x)f(x)d\sigma(x).\]
However on $S^{n-1}(y,r)$ the vector $n(x)=\frac{y-x}{\|x-y\|}$. So
$G(x-y)n(x)=\frac{1}{r^{n-1}}$. So
\[\int_{S^{n-1}(y,r)}G(x-y)n(x)f(x)d\sigma(x)=\int_{S^{n-1}(y,r)}
\frac{1}{r^{n-1}}(f(x)-f(y))d\sigma(x)\]
\[+\int_{S^{n-1}(y,r)}\frac{f(y)}{r^{n-1}}
d\sigma(x).\]
The right side of this previous expression reduces to
\[\int_{S^{n-1}(y,r)}\frac{(f(x)-f(y))}{r^{n-1}}d\sigma(x)+f(y)\int_{S^{n-1}}
d\sigma(x).\]
Now $\int_{S^{n-1}}d\sigma(x)=\omega_{n}$ and by continuity
$\lim_{r\rightarrow 0}\int_{S^{n-1}(y,r)}
\frac{(f(x)-f(y)}{r^{n-1}}d\sigma(x)=0$. The result follows. $\Box$

\ One important feature is to note that Kelvin inversion, ie $x^{-1}=\frac{-x}{\|x\|^{2}}$ whenever $x$ is non-zero, plays a
fundamental
role in this proof. Moreover the proof is almost exactly the same as the
proof of Cauchy's Integral Formula for piecewise $C^{1}$ curves in one
variable complex analysis.

\ Having obtained a Cauchy Integral Formula in ${\bf{R}}^{n}$ a number of
basic results that one might see in a first course in one variable complex
analysis carry over more or less automatically to the context described here.
This includes a Liouville Theorem and Weierstrass Convergence Theorem. We
leave it as an exercise to the interested reader to set up and establish the
Clifford analysis analogues of these results. Their statements and proofs can
be found in \cite{bds}.

\ Theorems 1 and 2 show us that the individual components of the equations
$Df=0$ and $gD=0$ comprise generalized Cauchy-Riemann equations. In the
particular case where $f$ is just vector valued so
$f=\Sigma_{j=1}^{n}f_{j}e_{j}$ then the generalized Cauchy-Riemann equations
become
$\frac{\partial f_{j}}{\partial x_{i}}=\frac{\partial f_{i}}{\partial x_{j}}$
whenever $i\neq j$ and $\Sigma_{j=1}^{n}\frac{\partial f_{j}}{\partial x_{j}}
=0$. This system of equations is often referred to as the Riesz system.

\ Having obtained an analogue of Cauchy's integral formula in euclidean
space we shall now exploit this result to show how many
consequences of the classical Cauchy integral carry over to the context
described here. We begin with the Mean Value Theorem.

\begin{theorem}
Suppose that $D(y, R)$ is a closed disc centered at $y$, of radius $R$ and
lying in $U$. Then for each left Clifford holomorphic function $f$ on $U$
\[f(y)=\frac{1}{R\omega_{n}}\int_{D(y, R)}\frac{f(x)}{\|x-y\|^{n-1}}dx^{n}.\]
\end{theorem}
{\bf{Proof:}}
We have already seen that for each $r\in(0, R)$
\[f(y)=\frac{1}{\omega_{n}}\int_{S^{n-1}(y, r)}\frac{f(x)}{\|x-y\|^{n-1}}
d\sigma(x),\]
where $S^{n-1}(y, r)$ is the $(n-1)$-dimensional sphere centered at $y$ and
of radius $r$. We obtain the result by integrating both sides of this
expression with respect to the variable $r$ and dividing throughout by $R$.
$\Box$

\ Let us now turn to explore the real analyticity properties of Clifford holomorphic
functions. First it may be noted that when $n$ is even
$G(x-y)=(-1)^{\frac{n-2}{2}}(x-y)^{-n+1}$. Also $(x-y)^{-1}=
x^{-1}(1-yx^{-1})^{-1}=(1-x^{-1}y)^{-1}x^{-1}$, and $\|x^{-1}y\|=\|yx^{-1}\|
=\frac{\|y\|}{\|x\|}$. So for $\|y\|<\|x\|$
\[(x-y)^{-1}=x^{-1}(1+yx^{-1}+\ldots+yx^{-1}\ldots yx^{-1}+\ldots)\]
\[=(1+x^{-1}y+\ldots+x^{-1}y\ldots x^{-1}y+\ldots)x^{-1}.\]
Hence these two sequences converge uniformly to $(x-y)^{-1}$ provided
$\|y\|\leq r<\|x\|$ and they converge pointwise to $(x-y)^{-1}$ provided
$\|y\|<\|x\|$. One can now take $(-1)^{\frac{n-2}{2}}$ times the $(n-1)$-fold
product of the series expansions of $(x-y)^{-1}$ with itself to obtain a
series expansion for $G(x-y)$. In this process of multiplying series
together in order to maintain the same radius of convergence one needs to
group together all linear combinations of monomials in $y_{1},\ldots, y_{n}$
that are of the same order. Thus we have deduced that when $n$ is even the
multiple Taylor series expansion
\[\Sigma_{j=0}^{\infty}(\Sigma_{\stackrel{j_{1},\ldots,j_{n}}{j_{1}+
\ldots+j_{n}=j}}
\frac{y_{1}^{j_{1}}\ldots y_{n}^{j_{n}}}{j_{1}!\ldots j_{n}!}
\frac{\partial^{j}G(x)}{\partial x_{1}^{j_{1}}\ldots \partial
x_{n}^{j_{n}}})\]
converges uniformly to $G(x-y)$ provided $\|y\|<r<\|x\|$ and converges
pointwise to $G(x-y)$ provided $\|y\|<\|x\|$.

\ A similar argument may be developed when $n$ is odd.

\ Returning to Cauchy's integral formula let us suppose that $f$ is a left
Clifford holomorphic function defined in a neighbourhood of the closure of some ball
$B(0, R)$. Then
\[f(y)=\frac{1}{\omega_{n}}\int_{\partial B(0,R)}G(x-y)n(x)f(x)d\sigma(x)=\]
\[\frac{1}{\omega_{n}}\int_{\partial B(0,R)}\Sigma_{j=0}^{\infty}
(\Sigma_{\stackrel{j_{1}\ldots j_{n}}{j_{1}+\ldots+j_{n}=j}}
\frac{y_{1}^{j_{1}}\ldots y_{n}^{j_{n}}}{j_{1}!\ldots j_{n}!}
\frac{\partial^{j}G(x)}{\partial x_{1}^{j_{1}}\ldots\partial x_{n}^{j_{n}}})
n(x)f(x)d\sigma(x)\]
provided $\|y\|<\|x\|$. As this series converges uniformly on each ball
$B(0, r)$ for each $r<R$ then this last integral can be re-written as
\[\frac{1}{\omega_{n}}\Sigma_{j=0}^{\infty}\int_{\partial B(0, R)}
(\Sigma_{\stackrel{j_{1}\ldots j_{n}}{j_{1}+\ldots+j_{n}=j}}
\frac{y_{1}^{j_{1}}\ldots y_{n}^{j_{n}}}{j_{1}!\ldots j_{n}!}
\frac{\partial^{j}G(x)}{\partial x_{1}^{j_{1}}\ldots \partial x_{n}^{j_{n}}}
n(x)f(x))d\sigma(x).\]
As the summation within the parentheses is a finite summation this last
expression easily reduces to
\[\frac{1}{\omega_{n}}\Sigma_{j=0}^{\infty}
(\Sigma_{\stackrel{j_{1}\ldots j_{n}}{j_{1}+\ldots+j_{n}=j}}
\frac{y_{1}^{j_{1}}\ldots y_{n}^{j_{n}}}{j_{1}!\ldots j_{n}!}
\int_{\partial B(0, R)}\frac{\partial^{j}G(x)}{\partial x_{1}^{j_{1}}\ldots
\partial x_{n}^{j_{n}}}n(x)f(x))d\sigma(x).\]
On placing
\[\frac{1}{\omega_{n}}\int_{\partial B(0, R)}\frac{\partial^{j}G(x)}
{\partial x_{1}^{j_{1}}\ldots\partial x_{n}^{j_{n}}}n(x)f(x)d\sigma(x)
=a_{j_{1}\ldots j_{n}}\]
it may be seen that on $B(0, R)$ the series
\[\Sigma_{j=0}^{\infty}
(\Sigma_{\stackrel{j_{1}\ldots j_{n}}{j_{1}+\ldots+j_{n}=j}}
\frac{x_{1}^{j_{1}}\ldots x_{n}^{j_{n}}}{j_{1}!\ldots j_{n}!}
a_{j_{1}\ldots j_{n}})\]
converges pointwise to $f(y)$. Convergence is uniform on each ball
$B(0, r)$ provided $r<R$.

\ Similarly if $g$ is a right Clifford holomorphic function defined in a neighbourhood
of the closure of $B(0, R)$ then the series
\[\Sigma_{j=0}^{\infty}
(\Sigma_{\stackrel{j_{1}\ldots j_{n}}{j_{1}+\ldots+j_{n}=j}}
b_{j_{1}\ldots j_{n}}\frac{y_{1}^{j_{1}}\ldots y_{n}^{j_{n}}}
{j_{1}!\ldots j_{n}!})\]
converges pointwise on $B(0, R)$ to $g(y)$ and converges uniformly on
$B(0, r)$ for $r<R$, where
\[b_{j_{1}\ldots j_{n}}=\frac{1}{\omega_{n}}\int_{\partial B(0, R)}
g(x)n(x)\frac{\partial^{j}G(x)}
{\partial x_{1}^{j_{1}}\ldots \partial x_{n}^{j_{n}}}d\sigma(x).\]

\ By translating the ball $B(0, R)$ to the ball $B(w, R)$ where
$w=w_{1}e_{1}+\ldots+w_{n}e_{n}$ one may readily observe that for any left
Clifford holomorphic function $f$ defined in a neighbourhood of the closure of $B(w, R)$
the series
\[\Sigma_{j=0}^{\infty}
(\Sigma_{\stackrel{j_{1}\ldots j_{n}}{j_{1}+\ldots j_{n}=j}}
\frac{(y_{1}-w_{1})^{j_{1}}\ldots (y_{n}-w_{n})^{j_{n}}}
{j_{1}!\ldots j_{n}!}a'_{j_{1}\ldots j_{n}})\]
converges pointwise on $B(w, R)$ to $f(y)$, where
\[a'_{j_{1}\ldots j_{n}}=\frac{1}{\omega_{n}}\int_{\partial B(w, R)}
\frac{\partial^{j}G(x-w)}{\partial x_{1}^{j_{1}}\ldots x_{n}^{j_{n}}}
n(x)f(x)d\sigma(x).\]
Again the series converges uniformly on $B(w, r)$ for each $r<R$.
A similar series may be readily obtained for any right Clifford holomorphic function
defined in a neighbourhood of the closure of $B(w, R)$.

\ The types of power series that we have developed for left Clifford holomorphic
functions are not entirely satisfactory. In particular, unlike their complex
analogues the homogeneous polynomials
\[\Sigma_{\stackrel{j_{1}\ldots j_{n}}{j_{1}+\ldots+j_{n}=j}}
\frac{x_{1}^{j_{1}}\ldots x_{n}^{j_{n}}}{j_{1}!\ldots j_{n}!}
a_{j_{1}\ldots j_{n}}\]
are not expressed as a linear combination of left Clifford holomorphic polynomials. To
rectify this situation let us first take a closer look at the Taylor
expansion for the Cauchy kernel $G(x-y)$ where all the Taylor coefficients are
real. Let us first look at the first order terms in the Taylor expansion.
This is the expression
\[y_{1}\frac{\partial G(x)}{\partial x_{1}}+\ldots+y_{n}
\frac{\partial G(x)}{\partial x_{n}}.\]
As $G$ is a Clifford holomorphic function then $\frac{\partial G(x)}{\partial x_{1}}=
-\Sigma_{j=2}^{n}e_{1}^{-1}e_{j}\frac{\partial G(x)}{\partial x_{j}}$.
Therefore the first order terms of the Taylor expansion for $G(x-y)$ can be
re-expressed as
\[\Sigma_{j=2}^{n}(y_{j}-e_{1}^{-1}e_{j}y_{1})\frac{\partial G(x)}
{\partial x_{j}}.\]
Moreover, for $2\leq j\leq n$ the first order polynomial
$y_{j}-e_{1}^{-1}e_{j}y_{1}$ is a left Clifford holomorphic polynomial.
Let us now go to second order terms. Again we will replace the operator
$\frac{\partial}{\partial x_{1}}$ by the operator $-\Sigma_{j=2}^{n}
e_{1}^{-1}e_{j}\frac{\partial}{\partial x_{j}}$ whenever it arises. Let us
consider the term $\frac{\partial^{2}G(x)}{\partial x_{i}\partial x_{j}}$
where $i\neq j\neq 1$. We end up with the polynomial
$y_{i}y_{j}-y_{i}y_{1}e_{1}^{-1}e_{j}-y_{j}y_{1}e_{1}^{-1}e_{i}
=\frac{1}{2}((y_{i}-y_{1}e_{1}^{-1}e_{i})(y_{j}-y_{1}e_{1}^{-1}e_{j})
+(y_{j}-y_{1}e_{1}^{-1}e_{j})(y_{i}-y_{1}e_{1}^{-1}e_{i})$.
Similarly the polynomial attached to the term $\frac{\partial^{2}G(x)}
{\partial x_{i}^{2}}$ is $(y_{i}-y_{1}e_{1}^{-1}e_{i})^{2}$. Using the
Clifford algebra anti-commutation relationship $e_{i}e_{j}+e_{j}e_{i}=-2
\delta_{ij}$ and on replacing the differential operator $\frac{\partial}
{\partial x_{1}}$ by the operator $-\Sigma_{j=2}^{n}e_{1}^{-1}e_{j}
\frac{\partial}{\partial x_{j}}$ it may be determined that the power
series we previously obtained for $G(x-y)$ can be replaced by the
series $\Sigma_{j=0}^{\infty}(\Sigma_{\stackrel{j_{2}\ldots j_{n}}
{j_{2}+\ldots+j_{n}=j}}P_{j_{2}\ldots j_{n}}(y)\frac{\partial^{j}G(x)}
{\partial x_{2}^{j_{2}}\ldots \partial x_{n}^{j_{n}}})$, where $\|y\|<
\|x\|$ and
\[P_{j_{2}\ldots j_{n}}(y)=\frac{1}{j!}\Sigma(y_{\sigma(1)}
-y_{1}e_{1}^{-1}e_{\sigma(1)})\ldots
(y_{\sigma(j)}-y_{1}e_{1}^{-1}e_{\sigma(j)}).\]
Here $\sigma(i)\in\{2,\ldots, n\}$ and the previous summation is taken over
all permutations of the
monomials $(y_{\sigma(i)}-y_{1}e_{1}^{-1}e_{\sigma(i)})$ without repetition.
The quaternionic monogenic analogues for these polynomials were introduced
by Fueter \cite{f} while the Clifford analogues, $P_{j_{2}\ldots j_{n}}$,
described here were introduced by Delanghe in \cite{d}. It should be noted that
each polynomial $P_{j_{2}\ldots j_{n}}(y)$ takes its values in the space
spanned by $1$, $e_{1}e_{2},\dots, e_{1}e_{n}$. Also each such polynomial is
homogeneous of degree $j$. Similar arguments to those
just outlined give that $G(x-y)=\Sigma_{j=0}^{\infty}
(\Sigma_{\stackrel{j_{2}\ldots j_{n}}{j_{2}+\ldots+j_{n}=j}}
\frac{\partial^{j}G(x)}{\partial x_{2}^{j_{2}}\ldots x_{n}^{j_{n}}}
\widetilde{P_{j_{2}\ldots j_{n}}}(y)$ provided $\|y\|<\|x\|$.
\begin{proposition}
Each of the polynomials $P_{j_{2}\ldots j_{n}}(y)$ is a left Clifford holomorphic
polynomial.
\end{proposition}
{\bf{Proof:}} As $DP_{j_{2}\ldots j_{n}}(y)=e_{1}(\frac{\partial}
{\partial y_{1}}+e_{1}^{-1}\Sigma_{j=2}^{n}e_{j}\frac{\partial}
{\partial y_{j}}P_{j_{2}\ldots j_{n}}(y))$ then we shall consider the
expression $(\frac{\partial}{\partial y_{1}}+\Sigma_{j=2}^{n}e_{1}^{-1}e_{j}
\frac{\partial}{\partial y_{j}})P_{j_{2}\ldots j_{n}}(y)$. This term is equal
to
\[(\frac{\partial}{\partial y_{j}}+\Sigma_{j=2}^{n}e_{1}^{-1}e_{j}
\frac{\partial}{\partial y_{j}})
\Sigma(y_{\sigma(1)}-e_{1}^{-1}e_{\sigma(1)}y_{1})\ldots\]
\[\ldots(y_{\sigma(i-1)}-e_{1}^{-1}e_{\sigma(i-1)}y_{1})
(y_{\sigma(i)}-e_{1}^{-1}
e_{\sigma(i)}y_{1})(y_{\sigma(i+1)}-e_{1}^{-1}e_{\sigma(i+1)}y_{1})\ldots\]
\[\ldots(y_{\sigma(j)}-e_{1}^{-1}e_{\sigma(j)}y_{1}).\]
This is equal to
\[\Sigma(y_{\sigma(1)}-e_{1}^{-1}e_{\sigma(1)}y_{1})\ldots
(y_{\sigma(i-1)}-e_{1}^{-1}e_{\sigma(i-1)}y_{1})(-e_{1}^{-1}e_{\sigma(i)})
(y_{\sigma(i+1)}-e_{1}^{-1}e_{\sigma(i+1)}y_{1})\]
\[\ldots (y_{\sigma(j)}-e_{1}^{-1}e_{\sigma(j)}y_{1})+
\Sigma e_{1}^{-1}e_{\sigma(i)}
(y_{\sigma(1)}-e_{1}^{-1}e_{\sigma(1)}y_{1})\ldots
(y_{\sigma(i-1)}-e_{1}^{-1}e_{\sigma(i-1)}y_{1})\]
\[(y_{\sigma(i+1)}-e_{1}^{-1}e_{\sigma(i+1)}y_{1})\ldots
(y_{\sigma(j)}-e_{1}^{-1}e_{\sigma(j)}y_{1}).\]

\ If we multiply the previous term by $y_{1}$ and add to it the following
term, which is equal to zero,
\[\Sigma(y_{\sigma(1)}-e_{1}^{-1}e_{\sigma(1)}y_{1})\ldots
(y_{\sigma(i-1)}-e_{1}^{-1}e_{\sigma(i-1)}y_{1})
(y_{\sigma(i)}-y_{\sigma(i)})\]
\[(y_{\sigma(i+1)}-e_{1}^{-1}e_{\sigma(i+1)}y_{1})\ldots
(y_{\sigma(j)}-e_{1}^{-1}e_{\sigma(i)}y_{1})\]
we get, after regrouping terms,
\[\Sigma(y_{\sigma(1)}-e_{1}^{-1}e_{\sigma(1)}y_{1})\ldots
(y_{\sigma(i-1)}-e_{1}^{-1}e_{\sigma(i-1)}y_{1})
(y_{\sigma(i)}-e_{1}^{-1}e_{\sigma(i)}y_{1})\]
\[(y_{\sigma(i+1)}-e_{1}^{-1}e_{\sigma(i+1)}y_{1})\ldots
(y_{\sigma(j)}-e_{1}^{-1}e_{\sigma(j)}y_{1})\]
\[-\Sigma(y_{\sigma(i)}-e_{1}^{-1}e_{\sigma(i)}y_{1})
(y_{\sigma(1)}-e_{1}^{-1}e_{\sigma(1)}y_{1})\ldots
(y_{\sigma(i-1)}-e_{1}^{-1}e_{\sigma(i-1)}y_{1})\]
\[(y_{\sigma(i+1)}-e_{1}^{-1}e_{\sigma(i+1)}y_{1})\ldots
(y_{\sigma(j)}-e_{1}^{-1}e_{\sigma(j)}y_{j}).\]
As summation is taken over all possible permutations without
repetition this last term vanishes. $\Box$

\ Using Proposition 1 and the results we previously obtained on series
expansions we can obtain the following generalization of Taylor expansions
from one variable complex analysis.
\begin{theorem}{\bf{(Taylor Series)}}
Suppose that $f$ is a left Clifford holomorphic function defined in an open
neighbourhood of the closure of the ball $B(w, R)$. Then
\[f(y)=\Sigma_{j=0}^{\infty}(\Sigma_{\stackrel{j_{2}\ldots j_{n}}
{j_{2}+\ldots+j_{n}=j}}P_{j_{2}\ldots j_{n}}(y-w)a_{j_{2}\ldots j_{n}}),\]
where $a_{j_{2}\ldots j_{n}}=\frac{1}{\omega_{n}}\int_{\partial B(w, R)}
\frac{\partial^{j}G(x-w)}{\partial x_{2}^{j_{2}}\ldots\partial x_{n}^{j_{n}}}
n(x)f(x)d\sigma(x)$ and $\|y-w\|<R$. Convergence is uniform
provided $\|x-w\|<r<R$.
\end{theorem}

\ A simple application of Cauchy's theorem now tells us that the Taylor
series that we obtained for $f$ in the previous theorem remains
valid on the largest open ball on which $f$ is defined and the largest open
ball on which $g$ is defined. Also the previous identities immediately yield
the mutual linear independence of the collection of the left Clifford holomorphic
polynomials $\{P_{j_{2}\ldots j_{n}}:j_{2}+\ldots j_{n}=j$ and $0\leq j
<\infty\}$.

\begin{center}{\Large Other Types of Clifford Holomorphic Functions}\end{center}

\ Unlike the the classical Cauchy-Riemann operator $\frac{\partial}
{\partial\overline{z}}=\frac{\partial}{\partial x}+i\frac{\partial}
{\partial y}$ the generalized Cauchy-Riemann operator $D$ that we have
introduced here does not have an identity component. Instead we could have
considered the differential operator
$D'=\frac{\partial}{\partial x_{0}}+\Sigma_{j+1}^{n-1}e_{j}\frac{\partial}
{\partial x_{j}}$. Also for $U'$ a domain in $R\oplus R^{n-1}$, the span of
$1, e_{1},\ldots, e_{n-1}$, one can consider $Cl_{n-1}$ valued differentiable
functions $f'$ and $g'$ defined on $U'$ such that $D'f'=0$ and $g'D'=0$,
where $g'D'=\frac{\partial g'}{\partial x_{0}}+\Sigma_{j=1}^{n-1}
\frac{\partial g'}{\partial x_{j}}e_{j}$. Traditionally such functions are
also called left monogenic and right monogenic functions. To avoid confusion
we shall call such functions unital left monogenic and unital right
monogenic respectively. In the case where $n=2$ the operator $D'$ corresponds
to the usual Cauchy-Riemann operator and unital monogenic functions are the
usual holomorphic functions studied in one variable complex analysis.
The function $G'(\underline{x})=\frac{\overline{\underline{x}}}
{\|\underline{x}\|^{n}}=\underline{x}^{-1}\|\underline{x}\|^{-n+2}$ is an
example of a function which is both unital
left monogenic and unital right monogenic. It is a simple matter to observe
that $f'$ is unital left monogenic if and only if $\tilde{f}'$ is unital
right monogenic. However $\overline{f}'$ is not unital right monogenic
whenever $f'$ is unital left monogenic. Instead $\overline{f}'$ satisfies
the equation $\overline{f}'D'=0$.

\ The function theory for unital left monogenic functions is much the same
as for left monogenic functions. For instance if $f'$ is unital left
monogenic on $U'$ and $g'$ is unital right monogenic on the same domain and
$S'$ is a piecewise smooth, compact surface lying in $U'$ and bounding a
subdomain $V'$ then $\int_{S'}g'(\underline{x})n(\underline{x})
f'(\underline{x})d\sigma(\underline{x})=0$ where $n(\underline{x})$ is the
outward pointing normal vector to $S'$ at $\underline{x}$. Also for each
$\underline{y}\in V'$ there is the following version of Cauchy's integral
formula
\[f'(\underline{y})=\frac{1}{\omega_{n}}\int_{S'}G'(\underline{x}-
\underline{y})n(\underline{x})f'(\underline{x})d\sigma(\underline{x}).\]

\ To get from the operator $D$ to the operator $D'$ one first rewrites $D$
as $e_{n}(\frac{\partial}{\partial x_{n}}+\Sigma_{j=1}^{n-1}e_{n}^{-1}e_{j}
\frac{\partial}{\partial x_{j}})$. On multiplying on the left by $e_{n}$ and
changing the variable $x_{n}$ to $x_{0}$ we get the operator
$D"=\frac{\partial}{\partial x_{0}}+\Sigma_{j=1}^{n-1}e_{n}^{-1}e_{j}
\frac{\partial}{\partial x_{j}}$. This operator takes its values in the
even subalgebra $Cl_{n}^{+}$ of $Cl_{n}$. On applying the isomorphism
\[\theta: Cl_{n-1}\rightarrow Cl_{n}^{+}:\theta(e_{j_{1}}\ldots e_{j_{r}})=e_{n}^{-1}e_{j_{1}}\ldots e_{n}^{-1}e_{j_{r}}\]
 it immediately follows that
$\theta(D')=D"$. So if $f'$ is unital left monogenic then $D"\theta(f)=0$.
If we change the variable $x_{0}$ of the function $\theta(f(\underline{x}))$
to $x_{n}$ we now get a left monogenic function, which we denote by
$\theta'(f)(x)$, where $x=x_{1}e_{1}+\ldots+x_{n}e_{n}\in U\subset R^{n}$
if and only if $\underline{x}=x_{n}+x_{1}e_{1}+\ldots+x_{n-1}e_{n-1}\in U'
\subset R\oplus R^{n-1}$.

\ It should be noted that $D'\overline{D'}=\overline{D'}D'=\triangle_{n}$.

\ When $n=3$ the algebra$Cl_{3}$ is split by the two projection
operators $E_{\pm}=\frac{1}{2}(1\pm e_{1}e_{2}e_{3})$ into a direct sum $E_{+}Cl_{3}\oplus E_{-}Cl_{3}$ and that
each of these subalgebras are isomorphic to the quaternion algebra ${\bf{H}}$.
In this setting the differential operator $E_{\pm}D'$ can best be written as
$\frac{\partial}{\partial t}+i\frac{\partial}{\partial x}+j\frac{\partial}
{\partial y}+k\frac{\partial}{\partial z}$ and the operator $E_{\pm}D$ can
best be written as $i\frac{\partial}{\partial x}+j\frac{\partial}
{\partial y}+k\frac{\partial}{\partial z}$. We shall denote the first of
these two operators by $D'_{{\bf{H}}}$ and the second by $D_{{\bf{H}}}$. The
operator $D'_{{\bf{H}}}$ is sometimes referred to as the
Cauchy-Riemann-Fueter operator. The function theory associated to the
differential operators $D_{{\bf{H}}}$ and $D'_{{\bf{H}}}$ is much the same as
that associated to the operators $D$ and $D'$. In fact historically the
starting point for Clifford analysis was to study the function theoretic
aspects of the operators $D'_{{\bf{H}}}$ and $D_{{\bf{H}}}$, see for instance
\cite{d,f} and the excellent review article of Sudbery \cite{su}.

\ Over all these function theories have proved itself to be much the same as
that for the operators $D$ and $D'$. It is a simple enough matter to set up
analogues of Cauchy's theorem and Cauchy's integral formula for the
quaternionic valued differentiable functions that are either annihilated by
$D'_{{\bf{H}}}$ or $D_{{\bf{H}}}$, either acting on the left or on the right.
When dealing with the operator $D'_{{\bf{H}}}$ such functions are called
quaternionic monogenic. The quaternionic monogenic Cauchy kernel is the
function $q^{-1}\|q\|^{-2}$. Consequently for each quaternionic left
monogenic function $f(q)$ defined on a domain $U"\subset{\bf{H}}$ and each
$q_{0}$ lying in a bounded subdomain with piecewise $C^{1}$ boundary  $S"$
\[f(q_{0})=\frac{1}{\omega_{3}}\int_{S"}(q-q_{0})^{-1}\|q-q_{0}\|^{-2}n(q)
f(q)d\sigma(q).\]
Similarly if $g$ is right quaternionic monogenic on $U"$ then
\begin{equation}
g(q_{0})=\frac{1}{\omega_{3}}\int_{S"}g(q)n(q)(q-q_{0})^{-1}\|q-q_{0}\|^{-2}
d\sigma(q).
\end{equation}

\begin{center}{\Large The Equation $D^{k}f=0$}\end{center}

\ It is reasonably well known that if $h$ is a real valued harmonic function
defined on a domain $U\subset R^{n}$ then for each $y\in U$ and each
compact, piecewise $C^{1}$ surface $S$ lying in $U$ such that $S$ bounds a
subdomain $V$ of $S$ and $y\in V$, then
\[h(y)=\frac{1}{\omega_{n}}\int_{S}(H(x-y)<n(x),grad h(x)>-<G(x-y),n(x)>h(x))
d\sigma(x),\]
where $H(x-y)=\frac{1}{(n-2)\|x-y\|^{n-2}}$. This formula is Green's formula
for a harmonic function, and it heavily relies on the standard inner product
on $R^{n}$. Introducing the Clifford algebra $Cl_{n}$ the right side  of
Green's formula is the real part of
\[\frac{1}{\omega_{n}}\int_{S}(G(x-y)n(x)h(x)-H(x-y)n(x)Dh(x))d\sigma(x).\]

\ Assuming that the function $h$ is $C^{2}$ then on applying Stokes' theorem
the previous integral becomes
\[\frac{1}{\omega_{n}}\int_{S^{n-1}(y, r(y))}
(G(x-y)n(x)h(x)-H(x-y)n(x)Dh(x))d\sigma(x),\]
where $S^{n-1}(y, r(y))$ is a sphere centered at $y$, of radius $r(y)$ and
lying in $V$. On letting the radius $r(y)$ tend to zero the first term of the
integral tends to $h(y)$ while the second term tends to zero. Consequently
the Clifford analysis version of Green's formula is
\[h(y)=\frac{1}{\omega_{n}}\int_{S}(G(x-y)n(x)h(x)-H(x-y)n(x)Dh(x))
d\sigma(x).\]

\ This formula was obtained under the assumption that $h$ is real valued and
$C^{2}$. The fact that we have assumed $h$ to be real valued can easily be
observed to be irrelevant, and so we can assume that $h$ is $Cl_{n}$ valued.
From now on we shall assume that all harmonic functions take their values in
$Cl_{n}$. If $h$ is also a left monogenic function then the Clifford analysis
version of Green's formula becomes Cauchy's integral formula.

\begin{proposition}
Suppose that $f$ is a Clifford holomorphic function on some domain $U$. Then $xf(x)$ is
harmonic.
\end{proposition}
{\bf{Proof:}}
$Dxf(x)=-nf(x)-\Sigma_{j=1}^{n}x_{j}\frac{\partial f(x)}{\partial x_{j}}
-\Sigma_{\stackrel{j, k}{j\neq k}}x_{k}e_{k}e_{j}
\frac{\partial f(x)}{\partial x_{j}}$.
Now
\[\Sigma_{\stackrel{j, k}{j\neq k}}x_{k}e_{k}e_{j}
\frac{\partial f(x)}{\partial x_{j}}=
\Sigma_{k=1}^{n}\Sigma_{j\neq k}x_{k}e_{k}e_{j}
\frac{\partial f(x)}{\partial x_{j}}.\]
As $f$ is left monogenic this last expression simplifies to
$\Sigma_{k=1}^{n}x_{k}\frac{\partial f(x)}{\partial x_{k}}$. Moreover
$D\Sigma_{j=1}^{n}x_{j}\frac{\partial f(x)}{\partial x_{j}}=0$. Consequently
$D^{2}xf(x)=0$. $\Box$

\ The previous proof is a generalization of the statement- "if $h(x)$ is a
real valued harmonic function then so is $<x,grad h(x)>$".

\ In fact in the previous proof we determine that
$Dxf(x)=-nf(x)-2\Sigma_{j=1}^{n}x_{j}\frac{\partial f(x)}{\partial x_{j}}$.
In the special case where $f(x)=P_{k}(x)$, a left Clifford holomorphic polynomial of
order $k$, this equation simplifies to $DxP_{k}(x)=-(n+2k)P_{k}(x)$. Suppose
now that $h(x)$ is a harmonic function defined in a neighbourhood of the
ball $B(0,R)$. Now $Dh$ is a left Clifford holomorphic function so we know that there
is a series $\Sigma_{l=0}^{\infty}P_{l}(x)$ of left Clifford holomorphic polynomials
with each $P_{l}$ homogeneous of degree $l$ and such that the series
converges locally uniformly on $B(0, R)$ to $Dh(x)$. Now consider
the series $\Sigma_{l=0}^{\infty}\frac{-1}{n+2l}P_{l}(x)$. As
$\frac{1}{n+2l}\|P_{l}(x)\|<\|P_{l}(x)\|$ then this new series converges
locally uniformly on $B(0, R)$ to a left Clifford holomorphic function $f_{1}(x)$.
Moreover, $Dxf_{1}(x)=Dh(x)$ on $B(0, R)$. Consequently $h(x)-xf_{1}(x)$ is
equal to a left Clifford holomorphic function $f_{2}(x)$ on $B(0, R)$. Thus we have
established:
\begin{proposition}
Suppose that $h$ is a harmonic function defined in a neighbourhood of
$B(0, R)$ then there are left Clifford holomorphic functions $f_{1}$ and $f_{2}$ defined
on $B(0, R)$ such that $h(x)=xf_{1}(x)+f_{2}(x)$ for each $x\in B(0, R)$.
\end{proposition}

\ This result remains invariant under translation. As a consequence
it shows us that all harmonic functions are real analytic functions. So
there is no need to specify whether or not a harmonic function is $C^{2}$.
The result also provides an Almansi type decomposition of harmonic functions
in terms of Clifford holomorphic functions over any ball in $R^{n}$.

\ It should be noted that Proposition 3 remains true if $h$ is only real
valued.

\ Proposition 3 gives rise to an alternative proof of the Mean Value
Theorem
for harmonic functions.

\begin{theorem}
For any harmonic function $h$ defined in a neighbourhood of a ball $B(a, R)$
\[h(a)=\frac{1}{\omega_{n}}\int_{\partial B(a, r)}h(x)d\sigma(x)\]
for any $r<R$.
\end{theorem}
{\bf{Proof:}}
 Proposition 3 tells us that there is a pair of left Clifford holomorphic functions
$f_{1}$ and $f_{2}$ such that $h(x)=(x-a)f_{1}(x)+f_{2}(x)$ on $B(a, R)$.
So $h(a)=f_{2}(a)$, and we have previously shown that $\frac{1}{\omega_{n}}
\int_{\partial B(a, r)}f_{2}(x)d\sigma(x)=f_{2}(a)$. Now
$\int_{\partial B(a,r)}(x-a)f_{1}(x)d\sigma(x)=
r\int_{\partial B(a, r)}n(x)f_{1}(x)d\sigma(x)=0$. $\Box$

\ The following is an immediate consequence of Proposition 3.
\begin{proposition}
If $h_{l}(x)$ is a harmonic polynomial homogeneous of degree $l$ then
\[h_{l}(x)=p_{l}(x)+xp_{l-1}(x)\]
where $p_{l}$ is a left Clifford holomorphic polynomial homogeneous of degree $l$ while
$p_{l-1}$ is a left monogenic polynomial which is homogeneous of degree $l-1$.
\end{proposition}

\ It is well known that pairs of homogeneous harmonic polynomials of
differing degrees of homogeneity are orthogonal with respect to the usual
inner product over the unit sphere. Proposition 4 offers a further
refinement
to this. Suppose that $f$ and $g$ are $Cl_{n}$ valued functions defined on
$S^{n-1}$ and each component of $f$ and $g$ is square integrable. If we define
the $Cl_{n}$ inner product of $f$ and $g$ to be
\[<f, g>=\frac{1}{\omega_{n}}\int_{S^{n-1}}\overline{f(x)}g(x)d\sigma(x)\]
then if $f$ and $g$ are both real valued this inner product is equal to
\[\frac{1}{\omega_{n}}\int_{S^{n-1}}f(x)g(x)d\sigma(x)\]
which is the usual inner product for real valued square integrable functions
defined on $S^{n-1}$. Now
\[<xp_{l-1}(x), p_{l}(x)>=
-\frac{1}{\omega_{n}}\int_{S^{n-1}}\overline{p}_{l-1}(x)xp_{l}(x)d\sigma(x)\]
\[=-\frac{1}{\omega_{n}}\int_{S^{n-1}}\overline{p}_{l-1}(x)n(x)p_{l}(x)
d\sigma(x)=0.\]
The evaluation of the last integral is an application of Cauchy's theorem.

\ Let us denote the space of $Cl_{n}$ valued functions defined on $S^{n-1}$
and such that each component is square integrable by $L^{2}(S^{n-1},Cl_{n})$.
Clearly the space of real valued square integrable functions defined on
$S^{n-1}$ is a subset of $L^{2}(S^{n-1},Cl_{n-1})$. The space
$L^{2}(S^{n-1},Cl_{n})$ is a $Cl_{n}$ module.

\ We have shown that by introducing the module $L^{2}(S^{n-1}, Cl_{n})$
Proposition 4 provides a further orthogonal decomposition of harmonic
polynomials using left Clifford holomorphic polynomials. We shall return to this theme
later. This decomposition was introduced for the case $n=4$ by Sudbery
\cite{su} and independently extended for all $n$ by Sommen \cite{s2}.

\ Let us now consider higher order iterates of the Dirac operator $D$. In the
same way as we have that $DH(x)=G(x)$ there is a function $G_{3}(x)$ defined
on $R^{n}\backslash\{0\}$ such that $DG_{3}(x)=H(x)$. Specifically
$G_{3}(x)=C(n, 3)\frac{x}{\|x\|^{n-2}}$ for some dimensional constant
$C(n,3)$. Continuing inductively we may find a function $G_{k}(x)$ on $R^{n}
\backslash\{0\}$ such that $DG_{k}(x)=G_{k-1}(x)$. Specifically
\[G_{k}(x)=C(n,k)\frac{x}{\|x\|^{n-k+1}}\]
when $n$ is odd and so is $k$.
\[G_{k}(x)=C(n,k)\frac{1}{\|x\|^{n-k}}\]
when $n$ is odd and $k$ is even
\[G_{k}(x)=C(n,k)\frac{x}{\|x\|^{n-k+1}}\]
when $n$ is even, $k$ is odd and $k<n$
\[G_{k}(x)=C(n,k)\frac{1}{\|x\|^{n-k}}\]
when $n$ is even, $k$ is even and $k<n$
\[G_{k}(x)=C(n,k)(x^{k-n}\log\|x\|+A(n,k)x^{k-n})\]
when $n$ is even and $k\geq n$. In the last expression $A(n,k)$ is a real
constant dependent on $n$ and $k$. $C(n,k)$ is a constant dependent on $n$
and $k$ throughout.

\ It should be noted that $G_{1}(x)=G(x)$ and $G_{2}(x)=H(x)$.  It should
also be noted that $D^{k}G_{k}(x)=0$.

\ Here is a simple technique for constructing solutions to the equation
$D^{k}g=0$ from left Clifford holomorphic functions. The special case $k=2$ was
illustrated in Proposition 2.
\begin{proposition}
Suppose that $f$ is a left Clifford holomorphic function on $U$ then
$D^{k}x^{k-1}f(x)=0$.
\end{proposition}
{\bf{Proof}} The proof is by induction. We have already seen the result to be
true in the case $k=2$ in Proposition 2. If $k$ is odd then $Dx^{k-1}f(x)
=(k-1)x^{k-2}f(x)$. If $k$ is even then
\[Dx^{k-1}f(x)=-n(k-1)x^{k-2}f(x)
+x^{k-2}\Sigma_{j=1}^{n}e_{j}x\frac{\partial f(x)}{\partial x_{j}}.\]
By arguments presented in Proposition 5 this expression is equal to
\[-n(k-1)x^{k-2}f(x)+x^{k-2}\Sigma_{j=1}^{n}x_{j}\frac{\partial f(x)}
{\partial x_{j}}.\]
The induction hypothesis tells us that the only term we need consider is
$x^{k-2}\Sigma_{j=1}^{n}x_{j}\frac{\partial f(x)}{\partial x_{j}}$.
However $\Sigma_{j=1}^{n}x_{j}\frac{\partial f(x)}{\partial x_{j}}$ is a left
Clifford holomorphic function. So proof by induction is now complete. $\Box$

\ In future we shall refer to a function $g:U\rightarrow Cl_{n}$ which
satisfies the equation $D^{k}g=0$ as a left $k$-monogenic function. Similarly
if $h:U\rightarrow Cl_{n}$ satisfies the equation $hD^{k}=0$ then $h$ is a
right $k$-monogenic function. In the case where $k=1$ we return to the
setting of left, or right, Clifford holomorphic functions and when $k=2$ we return to the
setting of harmonic functions. When $k=4$ the equations $D^{4}g=0$ and
$gD^{4}=0$ correspond to the equations $\triangle_{n}^{2}g=0$ and
$\triangle_{n}^{2}h=0$. So left or right $4$-monogenic functions are in fact
biharmonic functions. In greater generality if $k$ is even then a left or
right $k$-monogenic function $f$ automatically satisfies the equation
$\triangle_{n}^{\frac{k}{2}}f=0$.

\begin{proposition}
Suppose that $p$ is a left $k$-monogenic polynomial homogeneous of degree $q$
then there are left Clifford holomorphic polynomials $f_{0},\ldots, f_{k-1}$ such that
\[p(x)=f_{0}(x)+\ldots+x^{k-1}f_{k-1}(x)\]
and each polynomial $f_{j}$ is homogeneous of degree $q-j$ whenever $q-j\geq
0$ and is identically zero otherwise.
\end{proposition}
{\bf{Proof:}} The proof is via induction on $k$. The case $k=2$ is
established
immediately after the proof of Proposition 2. Let us now consider $Dp(x)$.
This is a left $k-1$-monogenic polynomial homogeneous of degree $q-1$. So by
the induction hypothesis $Dp(x)=g_{1}(x)+\ldots+x^{k-2}g_{k-1}(x)$ where
each $g_{j}$ is a left Clifford holomorphic polynomial homogeneous of degree $q-j$
whenever $q-j\geq 0$ and is equal to zero otherwise. Using Euler's lemma and
the observations made after the proof of Proposition 5 one may now find
left Clifford holomorphic polynomials $f_{1}(x),\ldots,f_{k-1}(x)$ such that
$D(xf_{1}(x)+\ldots+x^{k-1}f_{k-1}(x))=Dp(x)$ and $f_{j}(x)=c_{j}g_{j}(x)$
for some $c_{j}\in R$ and for $1\leq j\leq k-1$. It follows that $p(x)-
\Sigma_{j=1}^{k-1}x^{j}f_{j}(x)$ is a left Clifford holomorphic polynomial $f_{0}$
homogeneous of degree $q$. $\Box$

\ One may now use Proposition 6 and the arguments used to establish
Proposition 3 to deduce:
\begin{theorem}
Suppose that $f$ is a left $k$-monogenic function defined in a neighbourhood
of the ball $B(0, R)$ then there are left monogenic functions $f_{0},\ldots,
f_{k-1}$ defined on $B(0,R)$ such that
$f(x)=f_{0}(x)+\ldots+x^{k-1}f_{k-1}(x)$ on $B(0,R)$.
\end{theorem}

\ Theorem 6 establishes an Almansi decomposition for left $k$-monogenic
functions in terms of left Clifford holomorphic functions over any open ball. It also
follows from this theorem that each left $k$-monogenic function is a real
analytic function. It is also reasonably well known that if $h$ is a
biharmonic function defined in a neighbourhood of $B(0,R)$ then there are
harmonic functions $h_{1}$ and $h_{2}$ defined on $B(0, R)$ and such that
$h(x)=h_{1}(x)+\|x\|^{2}h_{2}(x)$. In the special case where $k=4$ Theorem
6 both establishes this result and refines it.

\ As each left $k$-monogenic function is a real analytic function then we
can immediately use Stokes' theorem to deduce the following Cauchy-Green
type formula.
\begin{theorem}
Suppose that $f$ is a left $k$-monogenic function defined on some domain
$U$ and suppose that $S$ is a piecewise $C^{1}$ compact surface lying in $U$
and bounding a bounded subdomain $V$ of $U$. Then for each $y\in V$
\[f(y)=\frac{1}{\omega_{n}}\int_{S}(\Sigma_{j=1}^{k}(-1)^{j-1}G_{j}(x-y)n(x)
D^{j-1}f(x))d\sigma(x).\]
\end{theorem}

\begin{center}{\Large Conformal Groups and Clifford Analysis}\end{center}

\ Here we will examine the role played
by the conformal group within parts of Clifford analysis. Our starting point
is to ask what type of diffeomorphisms acting on subdomains of $R^{n}$
preserve Clifford holomorphic functions. If a diffeomorphism $\phi$ can transform the
class of left Clifford holomorphic functions on one domain $U$ to a class of left
Clifford holomorphic functions on the domain $\phi(U)$ and do the same for the class
of right Clifford holomorphic functions on $U$ then it must preserve Cauchy's theorem.
So if $f$ is left Clifford holomorphic on $U$ and $g$ is right Clifford holomorphic on $U$ and
these functions are transformed to $f'$ and $g'$ respectively left and right
Clifford holomorphic functions on $\phi(U)$ then
\[\int_{S}g(x)n(x)f(x)d\sigma(x)=0=\int_{\phi(S)}g'(y)n(y)f'(y)d\sigma(y)\]
where $S$ is a piecewise $C^{1}$ compact surface lying in $U$ and $y=\phi(x)$.
An important point to note here is that we need to assume that $\phi$
preserves vectors orthogonal to the tangent spaces at $x$ and $\phi(x)$. As
the choice of $x$ and $S$ is arbitrary it follows that the diffeomorphism
$\phi$ is angle preserving. In other words $\phi$ is a conformal
transformation. A theorem of Liouville \cite{l} tells us that for dimensions
$3$ and greater the only conformal transformations on domains are
M\"{o}bius transformations.

\ In order to deal with M\"{o}bius transformations using Clifford algebras
we have seen in a previous chapter that one can use Vahlen matrices. We 
now proceed to show that each M\"{o}bius transformation
preserves monogenicity. Sudbery \cite{su} and also Bojarski \cite{b} have
established this fact. We will need the following two lemmata.

\begin{lemma}
Suppose that $\phi(x)=(ax+b)(cx+d)^{-1}$ is a M\"{o}bius transformation
then
\[G(u-v)=J(\phi, x)^{-1}G(x-y)\tilde{J}(\phi, y)^{-1}\]
where $u=\phi(x)$, $v=\phi(y)$ and $J(\phi, x)=\frac{(\widetilde{cx+d})}
{\|cx+d\|^{n}}$.
\end{lemma}
{\bf{Proof}} The proof essentially follows from the fact that
\[(x^{-1}-y^{-1})=x^{-1}(y-x)y^{-1}.\]
Consequently $\|x^{-1}-y^{-1}\|=\|x\|^{-1}\|x-y\|\|y\|^{-1}$.
Also $ax\tilde{a}-ay\tilde{a}=a(x-y)\tilde{a}$.

\ If one breaks the transformation down into terms arising from the generators
of the M\"{o}bius group and use the previous set of equations then one will
readily arrive at the result. $\Box$
\begin{lemma}
Suppose that $y=\phi(x)=(ax+b)(cx+d)^{-1}$ is a M\"{o}bius transformation
and for domains $U$ and $V$ we have $\phi(U)=V$ then
\[\int_{S}f(u)n(u)g(u)d\sigma(u)=\int_{\psi^{-1}(S)}f(\psi(x))\tilde{J}(\psi,x)
n(x)J(\psi,x)g(\psi(x))d\sigma(x)\]
where $u=\psi(x)$, $S$ is a orientable hypersurface lying in $U$ and
$J(\psi,x)=\frac{\widetilde{cx+d}}{\|cx+d\|^{n}}$.
\end{lemma}
{\bf{Outline Proof}} On breaking $\psi$ up into the generators of the
M\"{o}bius group the result follows from noting that
\[\frac{\partial x^{-1}}{\partial x_{j}}=-x^{-1}e_{j}x^{-1}.\]
$\Box$

\ It follows from Cauchy's Theorem that if $g(u)$ is a left Clifford holomorphic
function in the variable $u$ then $J(\psi,x)f(\psi(x))$ is left Clifford holomorphic
in the variable $x$.

\ When $\phi(x)$ is the Cayley transformation $y=(e_{n}x+1)(x+e_{n})^{-1}$
we can use this transformation to establish a Cauchy-Kowalewska extension
in a neighbourhood of the sphere. If $f(x)$ is a real analytic function
defined on an open subset $U$ of $S^{n-1}\backslash\{e_{n}\}$ then $l(y)=
J(\phi^{-1},y)^{-1}f(\phi(y))$ is a real analytic function on the open set
$V=\phi^{-1}(U)$. This function has a Cauchy-Kowalewska extension to a left
Clifford holomorphic function $L(y)$ defined on an open neighbourhood $V(g)\subset R^{n}$
of $V$. Consequently $F(x)=J(\phi^{-1}, x)L(\phi^{-1}(x))$ is a left
Clifford hholomorphic defined on an open neighbourhood $U(f)=\phi^{-1}(V(g))$ of $U$.
Moreover $F_{|U}=f$. Combing with similar arguments for the other Cayley
transformation $y=(-e_{n}x+1)(x-e_{n})^{-1}$ one can deduce:
\begin{theorem}
{\bf{(Cauchy-Kowalewska Theorem)}}
Suppose that $f$ is a $Cl_{n}$ valued real analytic function defined on
$S^{n-1}$. Then there is a unique left Clifford holomorphic function $F$ defined on an
open neighbourhood $U(f)$ of $S^{n-1}$ such that $F_{|S^{n-1}}=f$.
\end{theorem}

\ In fact if $f(u)$ is defined on some domain and satisfies the equation
$D^{k}f=0$ then the function $J_{k}(\psi,x)f(\psi(x))$ satisfies the same
equation, where $J_{k}(\psi,x)=\frac{\widetilde{cx+d}}{\|cx+d\|^{n-k+1}}$.
\begin{theorem}{\bf{(Fueter-Sce Theorem)}}
Suppose that $f=u+iv$ is a holomorphic function on a domain $\Omega\subset
{\bf{C}}$ and that $\Omega=\overline{\Omega}$ and $f(\overline{z})=
\overline{f(z)}$. Then the function $F(\underline{x})=u(x_{1},\|x'\|)
+e_{1}^{-1}\frac{x'}{\|x'\|}v(x_{1},\|x'\|)$ is a unital left
$n-1$-monogenic function
on the domain $\{\underline{x}:x_{1}+i\|x'\|\in\Omega\}$ whenever $n$ is
even. Here $x'=x_{2}e_{2}+\ldots+x_{n}e_{n}$.
\end{theorem}
{\bf{Proof:}}
First let us note that $x^{-1}e_{1}$ is left $n-1$ monogenic whenever
$n$ is even. It follows that $\frac{\partial^{k}}{\partial
x_{1}^{k}}x^{-1}e_{1}=c_{k}x^{-k-1}e_{1}$ is $n-1$ left monogenic for each
positive integer $k$. Here $c_{k}$ is some non-zero real number. Using
Kelvin inversion it follows that $x^{k}e_{1}$ is left $n-1$ monogenic for
each positive integer $k$. By taking translations and Taylor series
expansions for the function $f$ the result follows. $\Box$

\ This result was first established for the case $n=4$ by Fueter, \cite{f}, see
also Sudbery \cite{su}. It was extended to all even dimensions by Sce \cite{sc},
though the methods used do not make use of the conformal group. This result has been applied in \cite{q1,q2} to study various types of singular integral operators acting on $L^{p}$ spaces of Lipschitz perturbations of the sphere.

\begin{center}{\Large Conformally Flat Spin Manifolds}\end{center}

\ The invariance of monogenic functions under M\"{o}bius transformations described in the previous section makes use of a conformal weight factor $J(\psi,x)$. This invariance can be seen as an automorphic form invariance. This leads to a natural generalization of the concept of a Riemann surface to the euclidean setting. A manifold $M$ is said to be conformally flat if there is an atlas ${\cal A}$ of $M$ whose transition functions are M\"{0}bius transformations. For instance via the Cayley transformations $(e_{n+1}x+1)(x+e_{n+1})^{-1}$ and $(-e_{n+1}x+1)(x-e_{n+1})^{-1}$ one can see that the sphere $S^{n}\subset R^{n+1}$ is an example of a conformally flat manifold. Another way of constructing conformally flat manifolds is to take a simply connected domain $U$ of $R^{n}$ and consider a Kleinian subgroup $\Gamma$ of the M\"{o}bius group that acts discontinuously on $U$. Then the factorization $U\backslash \Gamma$ is a conformally flat manifold. For instance let $U=R^{n}$ and let $\Gamma$ be the integer lattice $Z^{k}=Ze_{1}+\ldots+Ze_{k}$ for some positive integer $k\leq n$. In this case $R^{n}\backslash Z^{k}$ gives the cylinder $C_{k}$ and when $k=n$ we get the $n$-torus. Also if we let $U=R^{n}\backslash\{0\}$ and $\Gamma=\{2^{k}:k\in Z\}$ the resulting manifod is $S^{1}\times S^{n-1}$. 

\ We locally construct a spinor bundle over $M$ by making the identification $(u,X)$ with either $(x,\pm J(\psi,x)X)$ where $u=\psi(x)=(ax+b)(cx+d)^{-1}=(-ax-b)(-cx-d)^{-1}$. If we can compatibly choose the signs then we have created a spinor bundle over the conformally flat manifold. Note, it might be possible to create more than one spinor bundle over $M$. For instance consider the cylinder $C_{k}$ if we make the identification $(x,X)$ with $x+\underline{m},(-1)^{m_{1}+\ldots+m_{l}}X)$ where $l$ is a fixed integer with $l\leq k$ and $\underline{m}=m_{1}e_{1}+\ldots+m_{l}e_{l}+\ldots+m_{k}e_{k}$ then we have created $k$ different spinor bundles $E^{1},\ldots E^{k}$ over $C_{k}$.

\ As we have used the conformal weight function function $J(\psi,x)$ to construct the spinor bundle $E$ then it is easy to see that a section $f:M\rightarrow E$ could be called a left monogenic section if it is locally a left monogenic function. It is now natural to ask if one can construct Cauchy integral formulas for such sections. To do this we essentially need to construct a kernel over the euclidean domain $U$ that is periodic with respect to $\Gamma$ and then use the projection map $p:U\rightarrow M$ to construct from this kernel a Cauchy kernel for $U$. In \cite{kr2} we show that the Cauchy kernel for  
$C_{k}$ with spinor bundle $E^{l}$ is constructed from the kernel
\[\cot_{k,l}(x,y)=\Sigma_{\underline{m}\in Z^{l},\underline{n}\in Z^{k-l}}(-1)^{m_{1}+\ldots+m_{l}}G(x-y+\underline{m}+\underline{n}),\]
where $\underline{n}=n_{l+1}e_{l+1}+\ldots+n_{n}e_{n}$.
While for the conformally flat spin manifold $S^{1}\times S^{n-1}$ with trivial bundle $Cl_{n}$ the Cauchy kernel is constructed from the kernel
\[\Sigma_{k=0}^{\infty}G(2^{k}x-2^{k}y)+2^{2-2n}G(x)(\Sigma_{k=-1}^{-\infty}G(2^{-k}x^{-1}-2^{-k}y^{-1}))G(y).\]

\ See \cite{k,kr1,kr2} for more details and related results. 

\ It should be noted that one may set up a Dirac operator over arbitrary Riemannian manifolds ,see for instance \cite{bbw}, and one may set up Cauchy integral formulas for functions annihilated by these Dirac operators, see for instance \cite{c,m2} for details. 

\begin{center}{\Large Boundary Behaviour and Hardy Spaces}\end{center}

\ Possibly the main topic that unites all that has been previously 
discussed here on Clifford analysis is its applications to boundary value 
problems. This in turn leads to a study of boundary behaviour of classes 
of Clifford holomorphic functions and Hardy spaces. Let us look first at one of the 
simplest cases. Previously we noted that if $\theta$ is a square integrable 
function defined on the sphere $S^{n-1}$ then there is a harmonic function 
$h$ defined on the unit ball in $R^{n}$ with boundary value $\theta$ 
almost everywhere. Previously we have seen that $h(x)=f_{1}(x)+xf_{2}(x)$ 
where $f_{1}$ and $f_{2}$ are left Clifford holomorphic. However on $S^{n-1}$ the 
function $G(x)$ equals $x$. So one can see that on $S^{n-1}$ we have 
$\theta(x)=f_{1}(x)+g(x)$ almost everywhere. Here $f_{1}$ is left 
monogenic on the unit ball $B(0,1)$ and $g$ is left Clifford holomorphic on 
$R^{n}\backslash clB(0,1)$, where $clB(0,1)$ is the closure of the open 
unit ball. Let $H^{2}(B(0,1))$ denote the space of Clifford holomorphic functions 
defined on $B(0,1)$ with extension to a square integrable function on 
$S^{n-1}$ and let $H^{2}(R^{n}\backslash clB(0,1)$ denote the class of 
left Clifford holomorphic functions defined on $R^{n}\backslash clB(0,1)$ with square 
integrable extension to the $S^{n-1}$. What we have so far outlined is 
that $L^{2}(S^{n-1})=H^{2}(B(0,1))\oplus H^{2}(R^{n}\backslash clB(0,1))$, 
where $L^{2}(S^{n-1})$ is the space of $Cl_{n}$ valued Lebesgue square 
integrable functions defined on $S^{n-1}$. This is the Hardy $2$-space 
decomposition of $L^{2}(S^{n-1})$. It is also true if we replace $2$ by 
$p$ where $1<p<\infty$. We will not go into this detail here as it is 
beyond the scope of the material covered here. 

\ Let us now take an alternative look at a way of obtaining this 
decomposition. This method will generalize to all reasonable surfaces. We 
will clarify what we mean by a reasonable surface later. Instead of 
considering an arbitrary square integrable function on $S^{n-1}$ let us 
instead assume that $\theta$ is a continuously differentiable function. Let 
us now consider the integral
\[\frac{1}{\omega_{n}}\int_{S^{n-1}}G(x-y)n(x)\theta(x)d\sigma(x)\]
where $y\in B(0,1)$. This defines a left Clifford holomorphic function on $B(0,1)$. 
Now let us allow the point $y$ to approach a point $z\in S^{n-1}$ along a 
differentiable path $y(t)$. Let us also assume that $\frac{dy(t)}{dt}$ 
evaluated at $t=1$, so $y(t)=z$, is not tangential to $S^{n-1}$ at $z$. We 
can essentially ignore this last point at a first read. We want to 
evaluate 
\[\lim_{t\rightarrow 
1}\frac{1}{\omega_{n}}\int_{S^{n-1}}G(x-y(t))n(x)\theta(x)d\sigma(x).\]
We do this by removing a small ball on $B(0,1)$ from $S^{n-1}$. The ball 
is centered at $z$ and is of radius $\epsilon$. We denote this ball by 
$b(z,\epsilon)$. The previous integral now splits into an integral over 
$b(z,\epsilon)$ and an integral over $S^{n-1}\backslash b(z,\epsilon)$. On 
$b(z, \epsilon)$ we can express $\theta(x)$ as 
$(\theta(x)-\theta(z))+\theta(z)$. As $\theta$ is continuously 
differentiable then $\|\theta(x)-\theta(z)\|<C\|x-z\|$ for some $C\in 
R^{+}$. It follows that 
\[\lim_{\epsilon\rightarrow 0}\lim_{t\rightarrow 
1}\int_{b(z,\epsilon)}\|G(x-y(t))n(x)(\theta(x)-\theta(z))\|d\sigma(x)=0.\]
Moreover, the term
\[\lim_{\epsilon\rightarrow 0}\lim_{t\rightarrow 
1}\frac{1}{\omega_{n}}\int_{b(z,\epsilon}G(x-y(t))n(x)\theta(z)d\sigma(x)\]
can, as $\theta(z)$ is a Clifford holomorphic function, be replaced by the term
\[\lim_{\epsilon\rightarrow 0}\lim_{t\rightarrow 
1}\int_{B(0,1)\cap\partial 
B(z,\epsilon)}G(x-y(t))n(x)\theta(z)d\sigma(x).\]
By the residue theorem the limit of this integral evaluates to 
$\frac{1}{2}\theta(z)$. 
 
\ We leave it to the interested reader to note that the singular integral 
or principal valued integral
\[\lim_{\epsilon\rightarrow 0}\lim_{t\rightarrow 
1}\frac{1}{\omega_{n}}\int_{S^{n-1}\backslash 
b(z,\epsilon)}G(x-y(t))n(x)\theta(x)d\sigma(x)=\]
\[P.V.\frac{1}{\omega_{n}}\int_{S^{n-1}}G(x-z)n(x)\theta(x)d\sigma(x)\]
is bounded. 

\ We have established that 
\[\lim_{t\rightarrow 1}\int_{S^{n-1}}G(x-y(t))n(x)\theta(x)d\sigma(x)=\]
\[\frac{1}{2}\theta(z)+
P.V.\frac{1}{\omega_{n}}\int_{S^{n-1}}G(x-z)n(x)\theta(x)d\sigma(x).\]

\ If we now assumed that $y(t)$ is a path tending to $z$ on the complement 
of $clB(0,1)$, then similar arguments give
\[\lim_{t\rightarrow 1}\int_{S^{n-1}}G(x-y(t))n(x)\theta(x)d\sigma(x)=\]
\[-\frac{1}{2}\theta(z)+
P.V.\frac{1}{\omega_{n}}\int_{S^{n-1}}G(x-z)n(x)\theta(x)d\sigma(x).\]
We will write these expressions as
\[(\pm\frac{1}{2}I+C_{S^{n-1}})\theta.\]

\ If we now consider the limit
\[\lim_{t\rightarrow 1}\frac{1}{\omega_{n}}\int_{S^{n-1}}G(x-y(t))n(x)
(\frac{1}{2}I+C_{S^{n-1}})\theta(x)d\sigma(x)\]
we may determine that 
\[(\frac{1}{2}I+C_{S^{n-1}})^{2}=\frac{1}{2}I+C_{S^{n-1}}.\]
Furthermore 
\[(\frac{1}{2}I+C_{S^{n-1}})(-\frac{1}{2}I+C_{S^{n-1}})=0\]
and
\[(-\frac{1}{2}I+C_{S^{n-1}})^{2}=-\frac{1}{2}I+C_{S^{n-1}}.\]

\ It is known that each function $\psi\in L^{2}(S^{n-1})$ can be 
approximated by a sequence of functions each with the same properties as 
$\theta$. This tells us that the previous formulas can be repeated but 
this time simply for $\theta\in L^{2}(S^{n-1})$. It follows that for such 
a $\theta$ we have
\[\theta=(\frac{1}{2}I+C_{S^{n-1}})\theta+(\frac{1}{2}I-C_{S^{n-1}})\theta.\]
This formula gives the Hardy space decomposition of $L^{2}(S^{n-1})$. In 
fact if one looks more carefully at the previous calculations to obtain 
these formulas we see that it is not so significant that the surface used 
is a sphere and we can re-do the calculations for any "reasonable" hypersurface 
$S$. In this case we get
\[\theta=(\frac{1}{2}I+C_{S})\theta+(\frac{1}{2}I-C_{S})\theta\]
where $\theta$ now belongs to $L^{2}(S)$ and 
\[C_{S}\theta=P.V.\frac{1}{\omega_{n}}\int_{S}G(x-y)n(x)\theta(x)d\sigma(x).\]
This gives rise to the Hardy space decomposition
\[L^{2}(S)=H^{2}(S^{+})\oplus H^{2}(S^{-})\]
where $S^{\pm}$ are the two domains that complement the surface $S$ (we 
are assuming that $S$ divides $R^{n}$ into two complementary domains.

\ Last one should address the smoothness of $S$. In some parts of the 
literature one simply assumes that $S$ is compact and $C^{2}$, while in 
more advanced and recent aspects of the literature one assumes that $S$ 
has rougher conditions, usually one assumes that the surface is simply 
Lipschitz continuous, see for instance \cite{lmq,lms,m1}. The formulas given above involving the singular 
integral operator $C_{S}$ are called Plemelj formulas. It is a simple 
exercise to see that these formulas are conformally invariant. So using 
Kelvin inversion or even a Cayley transformation one can see that these 
formulas and the Hardy space decompositions are also valid on unbounded 
surfaces and domains. A great deal of modern Clifford analysis has been 
devoted to the study of such Hardy spaces and singular integral operators. This is in fact due to an idea of R. Coifman that various hard problems in classical harmonic analysis studied in euclidean space might be more readily handled using tools from Clifford analysis, particularly the singular Cauchy transform and associated Hardy spaces. In particular Coifman speculated that a more direct proof of the celebrated Coifman-McIntosh-Meyer Theorem \cite{cmm}could be derived using Clifford analysis. The Coifman-McIntosh-Meyer Theorem establishes the $L^{2}$ boundedness of the double layer potential operator for Lipschitz graphs in $R^{n}$. Coifman's observation was that the double layer potential operator is the real or scalar part of the singular Cauchy transform arising in Clifford analysis and discussed earlier. If one can establish the $L^{2}$ boundedness of the singular Cauchy transform for a Lipschitz graph in $R^{n}$ then one automatically has the $L^{2}$ boundedness for the double layer potential operator for the same graph. The $L^{2}$ boundedness of the singular Cauchy transform was first established for Lipschitz graphs with small constant by M. Murray \cite{m} and extended to the general case by A. McIntosh, see for instance \cite{mc}, see also \cite{m1}. One very importand reason for needing to know that the double layer potential operator is $L^{2}$ bounded for Lipschitz graphs is to be able to solve bounndary value problems for domains with Lipschitz graphs as boundaries. such boundary value problems would include the Dirichlet problem and Neuaman problem for the Laplacian. See \cite{mc,m1} for more details. Moreover in \cite{w} Clifford analysis and more precisely the Hardy space decomposition mentioned here is specifically used to solve the water wave problem in three dimensions.

\begin{center}{\Large More on Clifford Analysis on the Sphere}\end{center}

\ In the previous section We saw that $L^{2}(S^{n-1})$ splits into a direct sum of Hardy spaces for the corresponding complemetary domains $B(0,1)$ and $R^{n}\backslash cl(B(0,1))$. In an earlier section we saw that any left Clifford holomorphic function $f(x)$ can be expressed as a locally uniformly convergent series $\Sigma_{j=0}^{\infty}f_{j}(x)$ where each $f_{j}(x)$ is left Clifford holomorphic and homogeneous of degree $j$. Now following \cite{su} consider the operator 
\[D=x^{-1}xD=x^{-1}(\Sigma_{i<k}e_{i}e_{k}(x_{i}\frac{\partial}{\partial x_{k}}-x_{k}\frac{\partial}{\partial x_{i}}-\Sigma_{j=1}^{n}x_{j}\frac{\partial}{\partial x_{j}}).\]
 By letting the last term in this expression act on homogeneous polynomials one may determine from Euler's lemma that in fact $\Sigma_{j=1}^{n}x_{j}\frac{\partial}{\partial x_{j}}$ is the radial operator $r\frac{\partial}{\partial r}$. So $r\frac{\partial}{\partial r}f_{j}(x)=jf_{j}(x)$. As each polynomial $f_{k}$ is Clifford holomorphic it follows that each $f_{j}$ is an eigenvector of the spherical Dirac operator $x\Lambda_{n-1}=x\Sigma_{i<k}e_{i}e_{k}(x_{i}\frac{\partial}{\partial x_{k}}-x_{k}\frac{\partial}{\partial x_{i}})$ with eigenvalue $k$. Now using Kelvin inversion we know that $f_{k}$ is  homogeneous of degree $k$ and left Clifford holomorphic if and only if $G(x)f_{k}(x^{-1})$ is homogeneous of degree $-n+1-k$ and is left Clifford holomorphic. On restricting $G(x)f_{k}(x^{-1})$ to the unit sphere this function becomes $xf_{k}(x^{-1})$ and this function is an eigenvector for the spherical Dirac operator $x\Lambda_{n-1}$. As each $f\in H^{2}(R^{n}\backslash cl(B(0,1)))$ can be written as $\Sigma_{k=0}^{\infty}G(x)f_{k}(x^{-1})$ where each $f_{k}$ is homogeneous of degree $k$ and is left Clifford holomorphic it follows that if $h\in L^{2}(S^{n-1})$ then $\Lambda_{n-1} xh(x)=(1-n)xh(x)-x\Lambda_{n-1} h(x)$. Similarly if we replace $S^{n-1}$ by the $n$-sphere $S^{n}$ embedded in $R^{n+1}$ then we have the identity $\Lambda_{n} xh(x)=-nxh(x)-x\Lambda_{n}h(x)$ for each $h\in L^{2}(S^{n})$. As all $C^{\infty}$ functions defined on $S^{n}$ belong to $L^{2}(S^{n})$ this identity holds for all such functions too.

\ It should be noted that for each $x\in S^{n}$ if we restrict the operator $x\Lambda_{n}$ to the tangent bundle $TS^{n}_{x}$ then we obtain the Euclidean Dirac operator acting on this tangent space.

\ By using the Cayley transformation $x=\psi(y)(e_{n+1}y+1)(y+e_{n+1})^{-1}$ from $R^{n}$ to $S^{n}\backslash\{e_{n+1}\}$ one can transform left Clifford holomorphic functions from domains in $R^{n}$ to functions defined on domains lying on the sphere. Namely if $f(y)$ is left Clifford holomorphic on the domain $U$ lying in $R^{n}$ then we obtain a function $f'(x)=J(\psi^{-1},x)f(\psi^{-1}(x)$ defined on the domain $U'=\psi(U)$ lying on $S^{n}$. Here $J(\psi^{-1},x)=\frac{x+1}{\|x+1\|^{n}}$. Similarly if $g(y)$ is right Clifford holomorphic on $U$ then $g'(x)=g(\psi^{-1}(x)J(\psi^{-1},x)$ is a well defined function on $U'$. Moreover for any smooth, compact hypersurface $S$ bounding a subdomain $V$ of $U$ we have from the conformal invariance of Cauchy's Theorem $\int_{S'}g'(x)n(x)f'(x)d\sigma'(x)=0$ where $S'=\psi(S)$, and $n(x)$ is the unit vector lying in the tangent space $TS^{n}_{x}$ of $S^{n}$ at $x$ and outer normal to 
$S'$ at $x$. Furthermore $sigma'$ is the Lebesgue measure on $S'$. 

\ From Lemma 1 it now follows that for each point $y'\in V'=\psi(V)$ we have the following version of Cauchy's Integral Formula:
\[f'(y')=\frac{1}{\omega_{n}}\int_{S'}G(x-y')n(x)f'(x)d\sigma(x)\]
where as before $G(x-y')=\frac{x-y'}{\|x-y'\|^{n}}$, but now $x$ and $y'\in S^{n}$. 
It would now appear that the functions $f'$ and $g'$ are solutions to some spherical Dirac equations. We need to isolate this Dirac operator. We shall achieve this by applying the operator $x\Lambda_{n}$ to the kernel $G_{s}(x,y')=G(x-y')$. As $x$ and $y'\in S^{n}$ then $\|x-y'\|^{2}=2-2<x,y'>$, where $<x,y'>$ is the inner product of $x$ and $y'$. So $G_{s}(x,y')=2^{\frac{-n}{2}}\frac{x-y'}{(1-<x,y'>)^{\frac{n}{2}}}$. So in calculating $x\Lambda_{n}G_{s}(x,y')$ we will need to know what $\Lambda_{n}<x,y'>$ evaluates to. It is a simple exercise to determine that $\Lambda_{n}<x,y'>=xy'+<x,y'>$ which is simply the wedge product, $x\wedge y'$, of $x$ with $y'$.  
Now let us calculate $x\Lambda_{n}G_{s}(x,y')$.
Now
\[x\Lambda_{n}G_{s}(x,y')=2^{\frac{-n}{2}}(x\Lambda_{n}\frac{x}{(1-<x,y'>)^{\frac{n}{2}}}-x\Lambda_{n}\frac{y'}{(1-<x,y'>)^{\frac{n}{2}}})\]
\[=2^{\frac{n}{2}}(-x\frac{nx}{(1-<x,y'>)^{\frac{n}{2}}}+\Lambda_{n}\frac{1}{(1-<x,y'>)^{\frac{-n}{2}}}-x\Lambda_{n}\frac{1}{(1-<x,y'>)^{\frac{n}{2}}}y').\]
First
\[\Lambda_{n}\frac{1}{(1-<x,y'>)^{\frac{n}{2}}}=\frac{n}{2}\frac{x\wedge y'}{(1-<x,y'>)^{\frac{n+2}{2}}}.\]
So
\[x\Lambda_{n}G_{s}(x,y')=2^{\frac{-n}{2}}\frac{n}{2(1-<x,y>)^{\frac{n+2}{2}}}(2(1-<x,y'>)+x\wedge y'-x(x\wedge y')y').\]
The expression 
\[2(1-<x,y'>)+x\wedge y'-x(x\wedge y')y'\]
is equal to
\[2-2<x,y'>+xy'+<x,y'>-x(xy+<x,y'>)y'.\]
This expression simplifies to
\[1-<x,y'>+xy'-xy'<x,y'>\]
which in turn simplifies to
\[(1-<x,y'>)(1+xy')=-x(1-<x,y'>)(x-y').\]
So
\[x\Lambda_{n}G_{s}(x,y')=\frac{-n}{2}xG_{s}(x,y').\]
Hence $x(\Lambda_{n}+\frac{n}{2})G_{s}(x,y')=0$. So the Dirac operator, $D_{s}$, over the sphere is $x(\Lambda_{n}+\frac{n}{2})$. It follows from our Cauchy integral formula for the sphere that $D_{s}f'(x)=0$. For more details on this operator, related operators and their properties see \cite{cnm,lr,r5,r6,v} 

\ Besides the operator $D_{s}$ we also need a Laplacian $\triangle_{s}$ acting on functions defined on domains on $S^{n}$. To do this we will work backwards, and look first for a fundamental solution to $\triangle_{n}$. A strong candidate for such a fundamental solution is the kernel $H_{s}(x,y')=\frac{1}{n-2}\frac{1}{\|x-y'\|^{n-2}}$. By similar considerations to those made in the previous calculation we find that $D_{s}H_{s}(x,y')=-xH_{s}(x,y')+G_{s}(x,y')$. So $(D_{s}+x)H_{s}(x,y')=G_{s}(x,y')$. Therefore we may define our Laplacian $\triangle_{s}$ to be $D_{s}(D_{s}+x)$.

\begin{definition}
Suppose $h$ is a $Cl_{n}$ valued function defined on a domain $U'$ of $S^{n}$. Then $h$ is called a harmonic function on $U'$ if $\triangle_{s}h=0$.
\end{definition}

\ In much the same way as one would derive Green's Theorem in $R^{n}$ one now has:
\begin{theorem}
Suppose $U'$ is a domain on $S^{n}$ and $h:U'\rightarrow Cl_{n}$ is a harmonic function on $U'$. Suppose also that $S'$ is a smooth hypersurface lying in $U'$ and that $S'$ bounds a subdomain $V'$ of $U'$ and that $y'\in V'$. Then
\[h(y')=\frac{1}{\omega_{n}}\int_{S'}(G_{s}(x,y')n(x)h(x)+H_{s}(x,y')n(x)D_{s}h(x))d\sigma'(x).\] 
\end{theorem}

\ See \cite{lr} for more details.

\begin{center}{\Large The Fourier Transform and Clifford Analysis}\end{center}

\ Closely related to Hardy spaces is the Fourier transform. Here we will consider $R^{n}$ as divided into upper and lower half space, $R^{n+}$ and $R^{n-}$, respectively. Where $R^{n+}=\{x=x_{1}e_{1}+\ldots+x_{n}e_{n}:x_{n}>o\}$ and $R^{n-}=\{X=x_{1}e_{1}+\ldots+x_{n}e_{n}:x_{n}<0\}$. These two domains have $R^{n-1}=$ span$<e_{1},\ldots,e_{n-1}>$ as common boundary. As before we have that $L^{2}(R^{n-1})=H^{2}(R^{n+})\otimes H^{2}(R^{n-})$. Let us now consider a function $\psi\in L^{2}(R^{n-1})$. Then \[\psi(y)=(\frac{1}{2}\psi(y)+\frac{1}{\omega_{n}}PV\int_{R^{n-1}}G(x'-y)e_{n}\psi(x')dx^{n-1})\]
\[+(\frac{1}{2}\psi(y)-\frac{1}{\omega_{n}}PV\int_{R^{n-1}}G(x'-y)e_{n}\psi(x')dx^{n-1})\]
almost everywhere. Here $\frac{1}{2}\psi(y)+\frac{1}{\omega_{n}}PV\int_{R^{n-1}}G(x'-y)e_{n}\psi(x')dx^{n-1}$ is the nontangential limit of $\frac{1}{\omega_{n}}\int_{R^{n-1}}G(x'-y(t))e_{n}\psi(x')dx^{n-1}$ as $y(t)$ tends to $y$ nontangentially through a smooth path in upper half space, and $\frac{1}{2}\psi(y)-\frac{1}{\omega_{n}}PV\int_{R^{n-1}}G(x'-y)e_{n}\psi(x')dx^{n-1}$ is the nontangential limit of $\frac{-1}{\omega_{n}}\int_{R^{n-1}}G(x'-y(t))e_{n}\psi(x')dx^{n-1}$ as $y(t)$ tends nontangentially to $y$ through a smooth path in lower half space. 

\ Consider now the Fourier transform, ${\cal F}(\psi)=\hat{\psi}$, of $\psi$. In order to proceed we need to calculate 
\[{\cal F}(\frac{1}{2}\psi\pm\frac{1}{\omega_{n}}PV\int_{R^{n-1}}G(x'-y)e_{n}\psi(x')dx^{n-1}).\]
In particular we need to determine ${\cal F}(\frac{1}{\omega_{n}}PV\int_{R^{n-1}}G(x'-y)e_{n}\psi(x')dx^{n-1})$. Following \cite{bla} it may be determined that this is $\frac{1}{2}i\frac{\zeta}{\|\zeta\|} e_{n}\hat{\psi}(\zeta)$ where $\zeta=\zeta_{1}e_{1}+\ldots+\zeta_{n-1}e_{n-1}$. So 
\[{\cal F}(\frac{1}{2}\psi\pm\frac{1}{\omega_{n}}\int_{R^{n-1}}G(x'-y)e_{n}\psi(x')d\sigma(x')=\frac{1}{2}(1\pm i\frac{\zeta}{\|\zeta\|} e_{n}).\]
Now as observed in \cite{blabla} 
\[(\frac{1}{2}(1\pm i\frac{\zeta}{\|\zeta\|} e_{n})^{2}=\frac{1}{2}(1\pm i\frac{\zeta}{\|\zeta\|} e_{n})\]
\[\frac{1}{2}(1\pm i\frac{\zeta}{\|\zeta\|} e_{n})\frac{1}{2}(1\mp i\frac{\zeta}{\|\zeta\|} e_{n})=0\]
and
\[i\zeta e_{n}\frac{1}{2}(1\pm i\frac{\zeta}{\|\zeta\|} e_{n})=\|\zeta\|\frac{1}{2}(1\pm i\frac{\zeta}{\|\zeta\|} e_{n}).\]
Now take the Cauchy Kowalewska extension of $e^{i<x',\zeta>}$. We get 
\[\exp(i<x',\zeta>-ix_{n}e_{n}\zeta)\]
defined on some neighbourhood in $R^{n}$ of $R^{n-1}$. Now consider
\[\exp(i<x',\zeta>-ix_{n}e_{n}\zeta)\frac{1}{2}(1\pm i\frac{\zeta}{\|\zeta\|} e_{n}).\]
This simplifies to 
\[e^{i<x',\zeta>+x_{n}\|\zeta\|}\frac{1}{2}(1+i\frac{\zeta}{\|\zeta\|} e_{n})\]
if $x_{n}>0$ and to
\[e^{i<x',\zeta>-x_{n}\|\zeta\|}\frac{1}{2}(1-i\frac{\zeta}{\|\zeta\|} e_{n})\]
if $x_{n}<0$. The first of these series converges locally uniformly on upper half space while the second series converges locally uniformly on lower half space. We denote these two functions by $e_{\pm}(x,\zeta)$ respectively. So the integrals $\frac{1}{\omega_{n}}\int_{R^{n-1}}e_{\pm}(x,\zeta)\hat{\psi}(\zeta)d\sigma(\zeta)$ define left monogenic functions $\Psi_{\pm}(x)$. These left monogenic functions are defined on upper and lower half space respectively. Moreover $\Psi_{\pm}\in H^{2}(R^{n\pm})$ respectively, and $\Psi_{\pm}$ explicitly give the Hardy space decomposition of $\psi\in L^{2}(R^{n-1})$. This approach to the links between the Fourier transform and Clifford analysis were first carried out in \cite{s3} and later independently rediscovered and applied in \cite{lmq}.

\ The integrals $\frac{1}{\omega_{n}}\int_{R^{n-1}}e_{\pm}(x,\zeta)d\zeta^{n-1}$ can be expressed in polar co-ordinates as an integral $\frac{1}{\omega_{n}}\int_{S^{n-2}}\int_{0}^{\infty}e^{ir<x',\zeta'>\pm x_{n}r}r^{n-2}\frac{1}{2}(1\pm\frac{\zeta}{\|\zeta\|}e_{n})drdS^{n-2}$, where $\zeta'=\frac{\zeta}{\|\zeta\|}$. In \cite{lmq} Chun Li observed that the integrals $\int_{0}^{\infty}e^{ir<x',\zeta'>\pm x_{n}r}\frac{1}{2}(1\pm \zeta' e_{n})r^{n-2}dr$ are Laplace transforms of the function $f(R)=R^{n-2}$. So the integral 
\[\int_{0}^{\infty}e^{ir<x',\zeta'>-x_{n}r}\frac{1}{2}(1-i\zeta' e_{n})dr\] 
evalutes to $(-i)^{n}(n-2)!(<x',\zeta'>+ix_{n})^{-n+1}$. Hence the integral 
\[\frac{1}{\omega_{n}}\int_{R^{n-1}}e_{+}(x,\zeta)d\zeta^{n-1}\] 
becomes 
\[\frac{1}{\omega_{n}}\int_{S^{n-2}}(-i)^{n}(n-2)!(<x',\zeta'>+ix_{n})^{-n+1}\frac{1}{2}(1+ie_{n}\zeta')dS^{n-2}.\]

\ For any complex number $a+ib$ the product $(a+ib)\frac{1}{2}(1+ie_{n}\zeta')$ is equal to $(a-be_{n})\frac{1}{2}(1+ie_{n}\zeta')$. Thus the previous integral becomes
\[\frac{1}{\omega_{n}}\int_{S^{n-2}}(-i)^{n}(<x',\zeta'>-x_{n}e_{n}\zeta')^{-n+1}\frac{1}{2}(1+ie_{n}\zeta')dS^{n-2}.\]
The imaginary, or $iCl_{n}$ part of this integral is the integral of an odd function so when $n$ is even the integral becomes
\[\frac{1}{\omega_{n}}\int_{S^{n-2}}\frac{(n-2)!e_{n}\zeta'}{2(<x',\zeta'>-x_{n}e_{n}\zeta')^{n-1}}dS^{n-2}\]
and when $n$ is odd the integral becomes
\[\frac{-1}{\omega_{n}}\int_{S^{n-2}}\frac{(n-2)!}{(<x',\zeta'>-x_{n}e_{n}\zeta')^{n-1}}dS^{n-2}.\]
These integrals are the plane wave decompositions of the Cauchy kernel for upper half space described by Sommen in \cite{s4}.
\ It should be noted that while introducing the Fourier transform and exploring some of its links with Clifford analysis we have also been forced to complexify the Clifford algebra $Cl_{n}$ so that we now work with the complex Clifford algebra $Cl_{n}(C)$. Further the functions $\frac{1}{2}(1\pm i\frac{\zeta}{\|\zeta\|e_{n}})$ can be seen as defined on spheres lying in the null cone $\{x_{n}e_{n}+iw': w'\in R^{n-1}$ and $x_{n}^{2}-\|w'\|^{2}=0\}$. This leads naturally to the question: what domains in $C^{n}$ do the functions $e_{\pm}(x,\zeta)$ extend to? Here we are replacing the real vector variable $x\in R^{n}$ by a complex vector variable $\underline{z}=z_{1}e_{1}+\ldots+z_{n}e_{n}\in C^{n}$, where $z_{1},\ldots,z_{n}\in C$. On placing $\underline{z}=x+iy$ where $x$ and $y$ are real vector variables, it may be noted that the term is $e^{-<\zeta,y'>-x_{n}\|\zeta\|}$ is well defined for $x_{n}\|\zeta\|>|<\zeta,y'>|$. So in this case $iy'+x_{n}e_{n}$ varies over the interior of the forward null cone $\{iy'+x_{n}e_{n}:x_{n}>0$ and $x_{n}>\|y'\|\}$. So $e_{+}(\underline{z},\zeta)$ is well defined for each $\underline{z}=x+iy=x'+iy'+(x_{n}+iy_{n})e_{n}\in C^{n}$ where $x'\in R^{n-1}$, $y_{n}\in R$, $x_{n}>0$ and $\|y'\|<x_{n}$. Similarly $e_{-}(x,\zeta)$ holomorphically extends to $\{\underline{z}=x'+iy'+(x_{n}+iy_{n})e_{n}:x'\in R^{n-1}$, $x_{n}<0$, $y_{n}\in R$ and $\|y'\|<|x_{n}|\}$. We denote these domains by $C^{\pm}$ respectively. It should be noted that $\Psi^{\pm}$ holomorphically continue to $C^{\pm}$ respectively. We denote these holomorphic continuations of $\Psi^{\pm}$ by $\Psi'^{\pm}$. The domains $C^{\pm}$ are examples of tube domains.

\begin{center}{\Large Complex Clifford Analysis}\end{center}

\ In the previous section we ended by showing that for any left Clifford holomorphic function $f\in H^{2}(R^{n,+})$ may be holomorphically continued to a function $f^{\dagger}$ defied on a tube domain $C^{+}$ in $C^{n}$. In this section we will briefly show how this type of holomorphic continuation happens for all Clifford holomorphic functions defined on a domain $U\subset R^{n}$. First let $S$ be a compact smooth hypersurface lying in $U$ and suppose that $S$ bounds a subdomain $V$ of $U$. Cauchy's integral formula gives us 
\[f(y)=\frac{1}{\omega_{n}}\int_{S}G(x-y)n(x)f(x)d\sigma(x)\]
for each $y\in V$. Let us now consider complexifying the Cauchy kernel. The function $G(x)$ holomorphically continues to $\frac{\underline{z}}{(\underline{z}^{2})^{\frac{n}{2}}}$. In even dimensions this is a well defined function on $C^{n}\backslash N(0)$ where $N(0)=\{\underline{z}\in C^{n}:\underline{z}^{2}=0\}$. In odd dimensions this lifts to a well defined function on complex $n$-dimensional Riemann surface double covering $C^{n}\backslash N(0)$. Though things work out well in odd dimensions we will for ease just work with the cases where $n$ is even. In holomorphically extending $G(x-y)$ in the variable $y$ we obtain a function $G^{\dagger}(x-\underline{z})=\frac{x-\underline{z}}{((x-\underline{z})^{2})^{\frac{n}{2}}}$. This function is well defined on $C^{n}\backslash N(x)$ where $N(x)=\{\underline{z}\in C^{n}:(x-\underline{z})^{2}=0\}$. So the integral $\frac{1}{\omega_{n}}\int_{S}G^{\dagger}(x-\underline{z})n(x)f(x)d\sigma(x)$ is well defined provided $\underline{z}$ is not in $N(x)$ for any $x\in S$. It may be determined that the set $C^{n}\backslash\cup_{x\in S}N(x)$ is an open set in $C^{n}$. We shall take the component of this open set which contains $V$. We shall denote this connected open subset of $C^{n}$ by $V^{\dagger}$. It now follows that the left Clifford holomorphic function $f(y)$ now has a holomorphic extension $f^{\dagger}(\underline{z})$ to $V^{\dagger}$. Furthermore this function is given by the integral formula 
\[f^{\dagger}(\underline{z})=\frac{1}{\omega_{n}}\int_{S}G^{\dagger}(x-\underline{z})n(x)f(x)d\sigma(x).\] 

\ Furthermore the holomorphic function $f^{\dagger}$ is now a solution to the complex Dirac equation $D^{\dagger}f^{\dagger}=0$, where $D^{\dagger}=\Sigma_{j=1}^{n}e_{j}\frac{\partial}{\partial z_{j}}$. By letting allowing the hypersurface to deform and move out to include more of $U$ in its interior we see that $f^{\dagger}$ is a well defined holomorphic function on $U^{\dagger}$ where $U^{\dagger}$ is the component of $C^{n}\backslash\cup_{x\in cl(U)\backslash U}N(x)$ which contains $U$. In the special cases where $U$ is either of $R^{n,\pm}$ then $U^{\dagger}=C^{\pm}$. See \cite{r1,r2,r3} for more details.

\end{document}